\documentclass{CSML}

\def\dOi{12(3:10)2016}
\lmcsheading%
{\dOi}
{1--13}
{}
{}
{Feb.~\phantom01, 2016}
{Sep.~14, 2016}
{}

\ACMCCS{[{\bf Theory of computation}]: Logic---Constructive mathematics; Models of computation---Computability---Recursive functions}


\usepackage{amssymb}
\usepackage{caption}
\usepackage{euler}
\usepackage{amsmath}
\usepackage{amsfonts}
\usepackage{graphicx}
\usepackage{tikz}
\usepackage{hyperref}

\usetikzlibrary{arrows}
\usetikzlibrary{calc}

\tikzset{citation/.style={font=\scriptsize\bfseries}}
\tikzset{impl/.style={-angle 45,rounded corners=1mm}}
\tikzset{biimpl/.style={angle 45-angle 45,rounded corners=1mm}}
\tikzset{plusimpl/.style={-angle 45,dashed}}
\tikzset{flowbox/.style={rectangle,draw=black,inner sep=2mm,align=center}}
\tikzset{equivclass/.style={flowbox,align=left,inner xsep=4mm,inner ysep=0mm}}

\makeatletter
\newcommand{\Hyp}[4][.5]{%
    \draw[impl] #2 -- #3
    let \p1 = ($ #3 - #2 $)
        in 
    \pgfextra{%
        \pgfmathparse{90-atan2(\x1,\y1)}
        \xdef\@ngle{\pgfmathresult}
    }
    node[citation,pos=#1,shift={(\@ngle+90:8pt)}]
              {\rotatebox{\@ngle}{#4}}
    }
\newcommand{\HyP}[4][.5]{%
    \draw[] #2 -- #3
    let \p1 = ($ #3 - #2 $)
        in 
    \pgfextra{%
        \pgfmathparse{90-atan2(\x1,\y1)}
        \xdef\@ngle{\pgfmathresult}
    }
    node[citation,pos=#1,shift={(\@ngle+90:8pt)}]
              {\rotatebox{\@ngle}{#4}}
    }
\makeatother

\newcommand{\equivclassbox}[2]{
\parbox{#1}{\begin{list}{}{
\setlength{\leftmargin}{0mm}
\setlength{\rightmargin}{0mm}
\setlength{\itemsep}{0mm}}
#2
\end{list}}
}

\theoremstyle{plain}
\newtheorem{theorem}{Theorem}

\newtheorem{lemma}[theorem]{Lemma}

\begin{document}

\title[Z-Stability in Constructive Analysis]{Z-Stability in Constructive Analysis\rsuper*}

\author[D.~Bridges]{Douglas Bridges}	
\address{School of Mathematics \& Statistics, University of Canterbury, Private Bag 4800, Christchurch 8140, New Zealand}	
\email{\{d.bridges, Maarten.Jordens\}@math.canterbury.ac.nz, semajdent@gmail.com}  

\author[J.~Dent]{James Dent}	
\address{\vspace{-18 pt}}	
\titlecomment{{\lsuper*}This research forms part of Dent's PhD thesis.}	

\author[M.~McKubre-Jordens]{Maarten McKubre-Jordens}	
\address{\vspace{-18 pt}}	


\maketitle

\begin{abstract}
We introduce Z-stability, a notion capturing the intuition that if a function $f$ maps a metric space into a normed space and if $\left\Vert f(x)\right\Vert $ is small, then $x$ is close to a zero of $f$. Working in Bishop's constructive setting, we first study pointwise versions of Z-stability and the related notion of good behaviour for functions. We then present a recursive counterexample to the classical argument for passing from pointwise Z-stability to a uniform version on compact metric spaces. In order to effect this passage constructively, we bring into play the positivity principle, equivalent to Brouwer's fan theorem for detachable bars, and the limited anti-Specker property, an intuitionistic counterpart to sequential compactness. The final section deals with connections between the limited anti-Specker property, positivity properties, and (potentially) Brouwer's fan theorem for detachable bars.
\end{abstract}

\section{Z-stability}

Let $f$ be a mapping of a metric space $\left(  X,\rho\right)  $ into a normed space $Y$, and let
\[ Z_{f}\equiv\left\{  x\in X:f(x)=0\right\}  =f^{-1}(0) \]
be the zero set of $f$. We say that $f$ is
\begin{itemize}
\item \emph{Z-stable at the point} $x\in X$ if for each $\varepsilon>0$ there exists $\delta>0$ such that if $\left\Vert f(x)\right\Vert <\delta$, then\footnote{{We adopt Richman's convention that, for example, }$\rho\left(  x,S\right)  <\varepsilon${\ means `there exists }$y\in S${\ with }$\rho\left(  x,y\right)  <\varepsilon${', and } $\rho\left(  x,S\right)  >0${\ means `there exists }$\delta >0${\ such that }$\rho\left(  x,y\right)  \geqslant\delta${\ for all }$y\in S${'.}} $\rho\left(  x,Z_{f}\right)  <\varepsilon$;
\item \emph{Z-stable }(on $X$) if it is Z-stable at each point of $X$; and
\item \emph{uniformly Z-stable on} $X$ if for each $\varepsilon>0$ there exists $\delta>0$ such that for each $x\in X$, if $\left\Vert f(x)\right\Vert<\delta$, then $\rho\left(  x,Z_{f}\right)  <\varepsilon$.
\end{itemize}

After preliminary work in the remainder of this introductory section, we move to a study of Z-stability and the related notion of good behaviour (well behavedness), defined in Section 2, within the framework of Bishop-style constructive mathematics (\textbf{BISH})---that is, roughly, mathematics with intuitionistic logic, an appropriate foundation such as CZF set theory \cite{Aczel} or Martin-L\"{o}f Type Theory \cite{ML}, and dependent choice.\footnote{For more on \textbf{BISH }see \cite{Bishop,BB,BR,BV}.} We then give a recursive counterexample to a classical theorem (see the next paragraph) connecting Z-stability and uniform Z-stability when $X$ is a compact metric space; this leads us to constructive counterparts of that classical theorem. In the final section we link the limited anti-Specker property (brought into play in the preceding section) and a pointwise positivity property for real-valued functions.

With classical logic it is trivial that $f$ is Z-stable on $X$: if $x\in X$, then either $f(x)=0$ and we may take $\delta=1$ in the Z-stability condition, or else $f(x)\neq0$ and we take $\delta=\left\Vert f(x)\right\Vert $. Moreover, if $f$ is continuous, then it is uniformly Z-stable on each compact subset of $X$. To see this, let $K\subset X$ be compact, let $\varepsilon>0$, and suppose that for each $\delta>0$ there exists $x\in K$ such that $\left\Vert f(x)\right\Vert <\delta$ and $\rho\left(  x,Z_{f}\right)\geqslant\varepsilon$. Then there exists a sequence $\left(  x_{n}\right)_{n\geqslant1}$ in $K$ such that for each $n$, $\left\Vert f(x_{n})\right\Vert<2^{-n}$ and $\rho\left(  x_{n},Z_{f}\right)  \geqslant\varepsilon$. Extracting a convergent subsequence, we may assume that $\left(  x_{n}\right)_{n\geqslant1}$ converges to a limit $x_{\infty}$ in $K$. Then $\rho\left(x_{\infty},Z_{f}\right)  \geqslant\varepsilon$; but $f$ is continuous, so $f(x_{\infty})=0$, a contradiction.

Things are not so straightforward in \textbf{BISH}: even cubic polynomial functions need not be Z-stable, as the following example shows. Let $\left(a_{n}\right)  _{n\geqslant1}$ be an increasing binary sequence, and define $f:\left(  -1,1\right)  \rightarrow\mathbf{R}$ by
\begin{equation}
f(x)\equiv x^{2}\left(  x-\frac{1}{2}\right)  -a=-a-\frac{1}{2}x^{2}+x^{3},
\label{aa1}
\end{equation}
where $a=\sum_{n=1}^{\infty}2^{-n}a_{n}$. First note that if $a>0$, then there is no zero of $f$ in the interval $\left[  -1/3,1/3\right]  $. Indeed, since $f^{\prime}(x)=3x^{2}-x=x\left(  3x-1\right)  $, we have $f^{\prime}(x)<0$, and therefore $f$ strictly decreasing, for $0<x<1/3$; for such $x$ it therefore follows that $f(x)<$ $f(0)=-a<0$; similarly, $f$ is strictly increasing on $[-1/3,0)$, and $f(x)<-a<0$ on that interval. Now suppose that $f$ is Z-stable at $0$. Compute $\delta>0$ such that if $\left\vert f(0)\right\vert <\delta$, then there exists $x$ such that $\left\vert x\right\vert <\frac{1}{3}$ and $f(x)=0$. Either $f(0)\neq0$, in which case $a\neq0$ and there exists $n$ with $a_{n}=1$, or else $\left\vert f(0)\right\vert <\delta$; in the latter case there exists $x\in\left(-1/3,1/3\right)  $ such that $f(x)=0$, which, by the foregoing, rules out the possibility that $a\neq0$ and therefore implies that $a=0$. We conclude that \emph{if every cubic function }$f:(-1,1)\rightarrow\mathbf{R}$\emph{\ is Z-stable, then we can derive} \emph{the essentially nonconstructive omniscience principle}

\begin{quote}
\textbf{LPO}: For every binary sequence $\left(  a_{n}\right)  _{n\geqslant1}$, either $a_{n}=0$ for all $n$ or else there exists (that is, we can compute) $N$ such that $a_{N}=1$.
\end{quote}

\noindent A variant of this example is used on page 3 of \cite{BV} to show how the standard classical interval-halving argument for root-finding with a computer will break down for a cubic function whose value at the midpoint of the interval under consideration is positive but smaller than the least positive number recognised by the computer. Given this, and that the intermediate value theorem in its full classical form is essentially nonconstructive, it is interesting to observe that Z-stability provides us with a constructively valid root-finding version of that theorem. Let $f:\left[  0,1\right]\rightarrow\mathbf{R}$ be Z-stable, sequentially continuous (we do not need full continuity in this argument), and such that $f(0)f(1)<0$. At each stage of the interval-halving argument, we can use the Z-stability of $f$ at the midpoint $m$ of the subinterval of $\left[  0,1\right]  $ under consideration, to decide either that there exists $y$ in that interval with $f(y)=0$ or else that $\left\vert f(m)\right\vert >0$; in the latter case, we proceed, as normally, to the next stage of the interval-halving.

These considerations illustrate why Z-stability might be something worthy of a constructive analysis:\footnote{{This is not to suggest that zero stability fails to have intrinsic merit as a constructive property.}} a common stopping criterion for root-finding algorithms is that the absolute value of the function be smaller than some predesignated positive quantity; but if the function is not known to be constructively Z-stable, then it is possible that its correct zeroes are not as close to the stopping point as we imagine. The same goes for algorithms for finding fixed points---such as that of Scarf for the Brouwer fixed-point theorem and, equivalently, the existence of equilibria in certain mathematical models of competitive economy \cite{Scarf,Scarf2}.

As so often happens, the particular problem caused by our cubic function---and hence, more generally, applicable to real-analytic functions---disappears when we enter the complex domain:
\begin{prop}
\label{0603a1}Let $K\subset\mathbf{C}$ be compact, with a totally bounded border $B$, and let $f:K\rightarrow\mathbf{C}$ a differentiable function such that $\inf_{z\in B}\left\vert f(z)\right\vert >\inf_{z\in K}\left\vert f(z)\right\vert $. Then $f$ is uniformly Z-stable on $K$.
\end{prop}
\proof
Corollary (5.3) of \cite[p. 153]{BB} tells us that $\inf_{z\in K}\left\vert f(z)\right\vert =0$; whence, by \cite[Theorem (5.11), p. 157]{BB}, there exist finitely many (not necessarily distinct) points $z_{1},\ldots,z_{m}$ of $K$, and a differentiable function $g:K\rightarrow\mathbf{C}$, such that
\[ f(z)=\left(  z-z_{1}\right)  \cdots(z-z_{m})g(z)\ \ \ \ \left(  z\in K\right) \]
and
\[ 0<\gamma\equiv\inf_{z\in K}\left\vert g(z)\right\vert . \]
Given $\varepsilon>0$, let $\delta\equiv\gamma\left(  \varepsilon/2\right)^{m}$. If $z\in K$ and $\left\vert f(z)\right\vert <\delta$, then
\[ \left\vert z-z_{1}\right\vert \cdots\left\vert z-z_{m}\right\vert <\left(\frac{\varepsilon}{2}\right)  ^{m} \]
and therefore there exists $k$ with $\left\vert z-z_{k}\right\vert<\varepsilon$.
\qed

From now on, we shall primarily be concerned with Z-stability in the more abstract context of metric and normed spaces.

\section{Z-stability and good behaviour}

Again let $f$ be a mapping of a metric space $X$ into a normed space $Y$. Consider the condition for Z-stability at a point $x\in X$. Written symbolically, that condition becomes
\[ \forall_{\varepsilon>0}\exists_{\delta>0}\left(  \left\Vert f(x)\right\Vert<\delta\Rightarrow\rho\left(  x,Z_{f}\right)  <\varepsilon\right)  . \]
We can rewrite this in the equivalent form\footnote{{This is not as obvious, constructively, as it may at first appear. The quantifiers are essential for the constructive equivalence to go through.}}
\[ \forall_{\varepsilon>0}\exists_{\delta>0}\left(  \rho\left(  x,Z_{f}\right)\geqslant\varepsilon\Rightarrow\left\Vert f(x)\right\Vert \geqslant\delta\right)  . \]
Classically, the $\varepsilon$ and $\delta$ are irrelevant; we might as well write
\[ \rho\left(  x,Z_{f}\right)  >0\Rightarrow f(x)\neq0. \]
Weakening this slightly, we obtain
\[ x\in\mathord{\sim}Z_{f}\Rightarrow f(x)\neq0. \]
where\footnote{{For elements }$x,y${\ of a metric space }$\left(  X,\rho\right)  ${, the expression `}$x\neq y${' stands for `}$\rho(x,y)>0${'. This is well known to be constructively stronger than `}$\lnot(x=y)$'.}
\[ \mathord{\sim}Z_{f}\equiv\left\{  y\in X:\forall_{z\in Z_{f}}\left(  y\neq z\right)\right\} \]
is the \emph{complement }of $Z_{f}$ in $X$. These considerations lead us to the definition: $f$ is \emph{well behaved} if $f(x)\neq0$ for each point $x\in\mathord{\sim}Z_{f}$.\footnote{{What about uniform Z-stability? The defining condition in that case can be rewritten as}
\[ \forall_{\varepsilon>0}\exists_{\delta>0}\forall_{x\in X}\left(  \rho\left(x,Z_{f}\right)  \geqslant\varepsilon\Rightarrow\left\Vert f(x)\right\Vert\geqslant\delta\right)  . \]
{This time, the }$\varepsilon${\ and }$\delta${\ are not irrelevant, even classically, and all we obtain is an alternative criterion for uniform Z-stability.}}

The proposition `every mapping of a metric space into a normed space is well behaved' is a constructive consequence of the statement
\[ \forall_{x\in\mathbf{R}}\left(  \lnot\left(  x=0\right)  \Rightarrow x\neq0\right)  . \]
which is equivalent to Markov's Principle,
\begin{quote}
\textbf{MP}: For each binary sequence $\left(  a_{n}\right)  _{n\geqslant1}$, if it is impossible that $a_{n}=0$ for all $n$, then there exists $N$ such that $a_{N}=1$,
\end{quote}
Since \textbf{MP} represents an unbounded search, we do not regard it as a valid principle of \textbf{BISH}. However, it is known that every linear mapping of a normed space \emph{onto} a Banach space is well behaved \cite[Theorem 1]{BI}, and that this proposition holds without the range being complete if and only if \textbf{MP} is derivable.

\begin{prop}
\label{0603a2}Let $f$ be a well-behaved mapping of a metric space $X$ into a normed space $Y$ such that
\[ Z_{f}\equiv\left\{  z\in X:f(z)=0\right\} \]
is inhabited and located in $X$. Then $f$ is Z-stable.
\end{prop}
\proof 
Consider a point $x\in X$ and $\varepsilon>0$. Either $\rho(x,Z_{f})<\varepsilon$ and we take $\delta=1$, or else $\rho\left(  x,Z_{f}\right)>0$. In the second case, since $f$ is well behaved, we have $f(x)\neq0$; taking $\delta=\left\Vert f(x)\right\Vert $, in view of \emph{ex falso quodlibet} we see that if $\left\Vert f(x)\right\Vert <\delta$, then $\rho\left(  x,Z_{f}\right)  <\varepsilon$.
\qed

Note that in the example of the cubic polynomial defined at (\ref{aa1}), the locatedness of $Z_{f}$ implies that either $a>0$ or $a=0$. For if $Z_{f}$ is located, then it is totally bounded, so
\[ \gamma\equiv\inf\left\{  x\in\left(  -1,1\right)  :f(x)=0\right\} \]
exists \cite[Corollary 2.2.7]{BV}. Either $\gamma>0$, in which case $a>0$, or else $\gamma<1/3$ and (as before) $a=0$.

For a converse to Proposition \ref{0603a2}, we introduce a completeness hypothesis.
\begin{prop}
\label{0603a3}A Z-stable mapping $f$ of a complete metric space $X$ into a normed space $Y$ is well behaved.
\end{prop}
\proof 
Fixing $x$ in $\mathord{\sim}Z_{f}$, construct a strictly decreasing sequence $\mathbf{\delta}$ converging to $0$ such that if $\left\Vert f(x)\right\Vert <\delta_{n}$, then $\rho\left(  x,Z_{f}\right)  <2^{-n}$. Then construct an increasing binary sequence $\mathbf{\lambda}$ such that
\begin{align*}
\lambda_{n}=0  &  \Rightarrow\left\Vert f(x)\right\Vert <\delta_{n},\\
\lambda_{n}=1  &  \Rightarrow\left\Vert f(x)\right\Vert >\delta_{n+1}\text{.}
\end{align*}
We may assume that $\lambda_{1}=0$. If $\lambda_{n}=0$, pick $x_{n}\in Z_{f}$ with $\rho\left(  x,x_{n}\right)  <2^{-n}$. If $\lambda_{n}=1-\lambda_{n-1}$, set $x_{k}=x_{n-1}$ for all $k\geqslant n$. Then $\left(  x_{n}\right)_{n\geqslant1}$ is a Cauchy sequence in $X$. To see this, consider $m,n$ with $m>n$. If $\lambda_{n}=1$, then $x_{m}=x_{n}$; if $\lambda_{m}=0$, then
\[ \rho\left(  x_{m},x_{n}\right)  \leqslant\rho\left(  x,x_{m}\right)+\rho\left(  x,x_{n}\right)  <2^{-m}+2^{-n}<2^{-n+1}. \]
It readily follows that $\rho\left(  x_{m},x_{n}\right)  <2^{-n+1}$ in all cases. Now, since $X$ is complete, the sequence $\left(  x_{n}\right)_{n\geqslant1}$ converges to a limit $x_{\infty}\in X$. Suppose that $f(x_{\infty})\neq0$. If there exists $n$ with $\lambda_{n}=1-\lambda_{n-1}$, then $x_{\infty}=x_{n-1}\in Z_{f}$, which is absurd. Hence $\lambda_{n}=0$ for all $n$, and therefore $f(x)=0$, which contradicts our original choice of $x$. We conclude that $\lnot\left(  f(x_{\infty})\neq0\right)  $, so $f(x_{\infty})=0$ and therefore $x\neq x_{\infty}$. Pick $N$ such that $\rho\left(x,x_{\infty}\right)  >2^{-N+1}$, and suppose that $\lambda_{N}=0$. If there exists $n>N$ such that $\lambda_{n}=1-\lambda_{n-1}$, then $x_{\infty}=x_{n-1}$, where $\rho\left(  x,x_{n-1}\right)  <2^{-n+1}\leqslant2^{-N+1}$. This contradicts our choice of $N$, so $\lambda_{n}=0$ for all $n\geqslant N$. Hence $x_{\infty}=\lim_{n\rightarrow\infty}x_{n}=x$, so $x\in Z_{f}$, a further contradiction. It follows that $\lambda_{N}\neq0$; whence $\lambda_{N}=1$ and therefore $\left\Vert f(x)\right\Vert >\delta_{N+1}$.
\qed

Examining the proof of Proposition \ref{0603a3} shows that in the hypotheses we can replace the completeness of $X$ by that of $Z_{f}$. Proposition \ref{0603a3} will be used in the proof of Proposition \ref{2303a2} in the next section.

\section{Uniform Z-stability}

We have already described a classical (in more than one sense) sequential compactness proof that a uniformly continuous mapping of a compact metric space into a normed space is uniformly Z-stable. We now present a recursive counterexample to that theorem. Since recursive constructive mathematics---that is, \textbf{BISH }supplemented by the Church-Markov-Turing thesis (CMT)---is (informally) a model of \textbf{BISH}, this example shows that in order to have a chance of deriving a good constructive counterpart of the classical theorem under consideration, we need to add to \textbf{BISH} some principles or hypotheses that run counter to CMT. That we shall do in due course.

Here, then, is our recursive example. Assuming CMT, we can construct a uniformly continuous, positive-valued function $g$ on $\left[  0,1/2\right]$ that has infimum $0$; see \cite{JR}, \cite[Ch. 6]{BR}, or \cite{BergerB2}. Let $f$ be the uniformly continuous mapping of $\left[  0,1\right]  $ into the nonnegative real line $\mathbf{R}^{0+}$ such that $f(x)=g(x)$ for $0\leqslant x\leqslant1/2$, $f(1)=0$, and $f$ is linear on $\left[1/2,1\right]  $. Then $Z_{f}=\left\{  1\right\}  $, which is inhabited and located. To prove that $f$ is Z-stable at each point $x$ of $\left[0,1\right]  $, let $\varepsilon>0$ and note that either $x<1$ or $x>1-\varepsilon$. In the first case, $f(x)>0$, so by \emph{ex falso}, if $\left\vert f(x)\right\vert <\frac{1}{2}\left\vert f(x)\right\vert $, then $\left\vert x-1\right\vert <\varepsilon$; this last inequality holds trivially in the second case. Now suppose that $f$ is uniformly Z-stable on $\left[0,1\right]  $. Then there exists $\delta>0$ such that if $x\in\left[0,1\right]  $ and $\left\vert f(x)\right\vert <\delta$, we can find $y\in Z_{f}$ such that $\left\vert x-y\right\vert <1/4$. But $\inf\left\{\left\vert f(x)\right\vert :0\leqslant x\leqslant1/2\right\}  =0$, so there exists $x\in\left[  0,1/2\right]  $ with $\left\vert f(x)\right\vert <\delta$. Clearly, $\rho\left(  x,Z_{f}\right)  =1-x>1/4$, a contradiction. This completes our recursive example of a uniformly continuous, Z-stable mapping $f:\left[  0,1\right]  \rightarrow\mathbf{R}$ (with inhabited, located zero set) that is not uniformly Z-stable.

This example hinges on the existence of a uniformly continuous, positive-valued function $f:\left[  0,1\right]  \rightarrow\mathbf{R}$ with infimum $0$. The natural addition to \textbf{BISH} that will counteract the recursive example is the case $X=\left[  0,1\right]  $ of the following \emph{positivity property}:
\begin{quote}
\textbf{POS}$^{X}$: \ If $f:X\rightarrow\mathbf{R}$ is uniformly continuous and positive-valued, then $\inf_{X}f$ (exists and) is positive.
\end{quote}
\begin{prop}
\label{2303a1} {The following are equivalent over \textbf{\emph{BISH}}}:
\begin{enumerate}[label=(\roman*)]
\item \textbf{\emph{POS}}$^{\left[  0,1\right]  }$.
\item \textbf{\emph{POS}}$^{2^{\mathbf{N}}}$.
\item \textbf{\emph{POS}}$^{X}$ holds for each compact metric
space $X$.
\end{enumerate}
\end{prop}
\proof 
Clearly, we need only prove that (i) $\Rightarrow$ (ii) and that (ii) $\Rightarrow$ (iii). First observe that, by a seminal result of Julian and Richman (\cite{JR}; see also \cite[Chapter 6]{BR}), Brouwer's fan theorem for detachable bars, \textbf{FT}$_{D}$, is equivalent to \textbf{POS}$^{\left[  0,1\right]  }$. Since (this is an easy exercise) \textbf{FT}$_{D}$ implies \textbf{POS}$^{2^{\mathbf{N}}}$, we see that (i) $\Rightarrow$ (ii). Next, assuming (ii) and given any compact metric space $X$, we apply Theorem (1.4) of \cite[Chapter 5]{BR}, to obtain a uniformly continuous mapping $g$ of $2^{\mathbf{N}}$ onto $X$. If $f:X\rightarrow\mathbf{R}$ is uniformly continuous and positive-valued, then so is $f\circ g:2^{\mathbf{N}}\rightarrow\mathbf{R}$; whence $\inf_{X}f=\inf_{2^{\mathbf{N}}}\left(  f\circ g\right)  >0$. Thus (ii) $\Rightarrow$ (iii).
\qed

We say that the implication (i) $\Rightarrow$ (iii) (respectively, (ii) $\Rightarrow$ (iii)) in Proposition \ref{2303a1} shows that \textbf{POS}$^{\left[  0,1\right]  }$ (respectively, \textbf{POS}$^{2^{\mathbf{N}}}$) is \emph{prototypical} for the positivity property on compact (metric) spaces.

This brings us to our first result on the passage from Z-stability to uniform Z-stability.
\begin{prop}
\label{2303a2}Let $X$ be a compact metric space with the positivity property, and $f$ a Z-stable, uniformly continuous mapping of $X$ into a normed space $Y$ such that $Z_{f}$ is inhabited and located. Then $f$ is uniformly Z-stable on $X$.
\end{prop}
\proof 
Given $\varepsilon>0$, and referring to \cite[Chapter 4, Theorem (4.9)]{BB}, we may assume without loss of generality that
\[ K\equiv\left\{  x\in X:\rho\left(  x,Z_{f}\right)  \geqslant\frac{\varepsilon}{2}\right\} \]
is compact. By Proposition \ref{0603a3}, since $f$ is Z-stable and $K$ is complete, the mapping $x\rightsquigarrow\left\Vert f(x)\right\Vert $ is positive-valued on $K$. Hence
\[ 0<\delta\equiv\inf_{x\in K}\left\Vert f(x)\right\Vert . \]
If $x\in X$ and $\left\Vert f(x)\right\Vert <\delta$, then $x\notin K$ and therefore $\rho\left(  x,Z_{f}\right)  <\varepsilon$.
\qed

\begin{cor}
\label{2303a3}\textbf{\emph{BISH + POS}}$^{\left[  0,1\right]  }\vdash$ Let $X$ be a compact metric space, and $f$ a Z-stable, uniformly continuous mapping of $X$ into a normed space $Y$ such that $Z_{f}$ is inhabited and located. Then $f$ is uniformly Z-stable on $X$.
\end{cor}
\proof 
This follows from Propositions \ref{2303a1} and \ref{2303a2}.
\qed

\section{Anti-Specker properties and Z-stability}

Now let $X$ be a subspace of a metric space $\left(  E,\rho\right)  $, and $\omega$ a point of $E$ with $\rho\left(  \omega,X\right)  >0$. We call the metric space $X\cup\left\{  \omega\right\}  $ a \emph{one-point extension} of $X$. It is straightforward to construct one-point extensions of a given metric space $X$.

Recall that a sequence $\left(  x_{n}\right)  _{n\geqslant1}$ in $E$ is said to be \emph{eventually bounded away from the point} $x\in E$ if there exist $N$ and $\delta>0$ such that $\rho(x,x_{n})\geqslant\delta$ for all $n\geqslant N$. Specker's theorem from recursive constructive analysis (see \cite[Chapter 3]{BR}) says that there exists a sequence in $\left[0,1\right]  $ that is eventually bounded away from each point of that interval; this is, of course, a strong recursive counterexample to the sequential compactness of $\left[  0,1\right]  $. Various antitheses of Specker's theorem have been studied as constructive substitutes for sequential compactness; see \cite{BergerB,dsb1,BD1,Dien,Dien1b,Dien13}. One of the weakest of those notions is the \emph{limited anti-Specker property} (relative to one-point extensions),
\begin{quote}
\textbf{AS}$_{L}^{X}$: If $X\cup\left\{  \omega\right\}  $ is a one-point extension of $X$, and $\left(  x_{n}\right)  _{n\geqslant1}$ is a sequence in $X\cup\left\{  \omega\right\}  $ that is eventually bounded away from each point of $X$, then there exists $n$ such that $x_{n}=\omega$.
\end{quote}
The property \textbf{AS}$_{L}^{X}$ is independent of the one-point extension of $X$ (cf. \cite[Propositions 1 and 2]{dsbarch}). It was introduced in \cite{Dent} and further discussed in \cite{BDMM}.

In Proposition \ref{2303a2} we can obtain the same conclusion if we replace the positivity property with the limited anti-Specker property and add separability:

\begin{prop}
\label{2003a1}Let $X$ be a compact metric space with the limited anti-Specker property, and $f$ a Z-stable, uniformly continuous mapping of $X$ into a normed space $Y$ such that $Z_{f}$ is separable. Then $f$ is uniformly Z-stable on $X$.
\end{prop}
\proof 
Let $X\cup\left\{  \omega\right\}  $ be a one-point extension of $X$ with $\rho\left(  \omega,X\right)  \geqslant1$, and fix $\varepsilon\in\left(0,1\right)  $. By \cite[Chapter 4, Theorem (4.9)]{BB}, there exists a strictly decreasing sequence $\mathbf{\delta}\equiv\left(  \delta_{n}\right)_{n\geqslant1}$ such that for each $n$, $0<\delta_{n}<2^{-n}$ and the set
\[ K_{n}\equiv\left\{  x\in X:\left\Vert f(x)\right\Vert \leqslant\delta_{n}\right\} \]
is compact. Let $\left(  z_{n}\right)  _{n\geqslant1}$ be a dense sequence in $Z_{f}$, and for each $n$ write
\[ S_{n}\equiv\left\{  z_{1},\ldots,z_{n}\right\} \]
Construct a binary sequence $\mathbf{\lambda}\equiv\left(  \lambda_{n}\right)_{n\geqslant1}$ such that
\begin{align*}
\lambda_{n}=0  &  \Rightarrow\sup\left\{  \rho\left(  x,S_{n}\right)  :x\in K_{n}\right\}  >\frac{\varepsilon}{2},\\
\lambda_{n}=1  &  \Rightarrow\sup\left\{  \rho\left(  x,S_{n}\right)  :x\in K_{n}\right\}  <\varepsilon.
\end{align*}
If $\lambda_{n}=0$, pick $x_{n}\in K_{n}$ with $\rho\left(  x_{n},S_{n}\right)  >\varepsilon/2$. If $\lambda_{n}=1$, set $x_{n}\equiv\omega$. We show that $\left(  x_{n}\right)  _{n\geqslant1}$ is eventually bounded away from each point $x$ of $X$. It will suffice to show that there exists $c>0$ such that $\rho\left(  x_{n},x\right)  \geqslant c$ for all sufficiently large $n$ with $\lambda_{n}=0$. By the Z-stability of $f$ at $x$, we can find $\alpha>0$ such that if $\left\Vert f(x)\right\Vert <\alpha$, then $\rho\left(  z,Z_{f}\right)  <\varepsilon/4$. Either $f(x)\neq0$ or $\left\Vert f(x)\right\Vert <\alpha$. In the first case, choose $\nu$ such that $\delta_{\nu}<\frac{1}{2}\left\Vert f(x)\right\Vert $, and then $\gamma\in\left(  0,1\right)  $ such that $\left\Vert f(x)-f(y)\right\Vert<\frac{1}{2}\left\Vert f(x)\right\Vert $ whenever $y\in X$ and $\rho\left(x,y\right)  <\gamma$. Then for all $n\geqslant\nu$ with $\lambda_{n}=0$ we have
\begin{align*}
\left\Vert f(x)-f(x_{n})\right\Vert  &  \geqslant\left\Vert f(x)\right\Vert -\delta_{n}\geqslant\left\Vert f(x)\right\Vert -\delta_{\nu}\\
&  >\left\Vert f(x)\right\Vert -\frac{1}{2}\left\Vert f(x)\right\Vert =\frac{1}{2}\left\Vert f(x)\right\Vert ,
\end{align*}
so $\rho\left(  x,x_{n}\right)  \geqslant\gamma$. In the case $\left\Vert f(x)\right\Vert <\alpha$, pick $\zeta\in Z_{f}$ with $\rho\left(\zeta,x\right)  <\varepsilon/4$; then choose $N$ such that $\rho\left(\zeta,z_{N}\right)  <\varepsilon/4-\rho\left(  \zeta,x\right)  $. For all $n\geqslant N$ with $\lambda_{n}=0$ we have $\rho(x_{n},z_{N})\geqslant\rho\left(  x_{n},S_{n}\right)  >\varepsilon/2$, so
\[ \rho\left(  x_{n},x\right)  \geqslant\rho\left(  x_{n},z_{N}\right)-\rho\left(  z_{N},\zeta\right)  -\rho\left(  \zeta,x\right)  >\frac{\varepsilon}{2}-\frac{\varepsilon}{4}=\frac{\varepsilon}{4}. \]
This completes the proof that $\left(  x_{n}\right)  _{n\geqslant1}$ is eventually bounded away from each point of $X$. We now apply the limited anti-Specker property in $X$, to compute $N$ with $x_{N}=\omega$ and therefore $\lambda_{N}=1$. For each $x\in X$ with $\left\Vert f(x)\right\Vert<\delta_{N}$, we have $x\in K_{N}$, so (as $\lambda_{N}=1$)
\[ \rho\left(  x,Z_{f}\right)  \leqslant\rho\left(  x,S_{N}\right)<\varepsilon. \]
Since $\varepsilon$ is arbitrary, we have shown that $f$ is uniformly Z-stable on $K$.
\qed

Next, with the aid of a stronger property than \textbf{AS}$_{L}^{X}$, we head towards Proposition \ref{aa3}, a generalised form of the principle of isolation of zeroes for complex analytic functions (cf. \cite[pp.194--195]{Dieudonne}). Although that proposition is not about uniform Z-stability, it uses both Z-stability and an anti-Specker property, and so is a fitting digression from the main theme of the section.

Let $\left(  X,\rho\right)  $ be a metric space. Recall that if $S,T$ are subsets of $X$ for which there exists $r>0$ such that if $x\in X$ and $\rho\left(  x,S\right)  <r$ entails $x\in T$, then $S$ is said to be \emph{well contained }in $T$, and we write $S\subset\subset T$. Let $U$ be an open subset of $X$, and $f$ a mapping of $U$ into a normed space $Y$. We say that the set $Z_{f}\subset U$ of zeroes of $f$ is \emph{countably isolated }if there is a one-one enumeration $\left(  z_{n}\right)  _{n\geqslant1}$ of $Z_{f}$ that is eventually bounded away from each of its terms $z_{m}$. In that case, any one-one enumeration of $Z_{f}$ is eventually bounded away from each of its terms.

\begin{lemma}
\label{aa2}Let $E$ be a complete metric space, $U$ an open subset of $E$, and $f$ a pointwise continuous mapping of $U$ into a normed space $Y$. Suppose that $f$ is Z-stable at each point of $U$, and that $Z_{f}$ is countably isolated, with one-one enumeration $z_{1},z_{2},\ldots$ . Then $\left(z_{n}\right)  _{n\geqslant1}$ is eventually bounded away from each point of $U$.
\end{lemma}
\proof 
Fixing $x\in U$, pick $r>0$ such that the closed ball $\overline{B}(x,r)$ is well contained in $U$. For each positive integer $n$ choose $\delta_{n}>0$ such that if $\left\Vert f(x)\right\Vert <\delta_{n}$, then there exists $\zeta$ with $f(\zeta)=0$ and $\rho\left(  \zeta,x\right)  <2^{-n}r$; we may assume that $\delta_{n+1}<\delta_{n}$. Construct an increasing binary sequence $\mathbf{\lambda}$ such that for each $k\geqslant1$,
\begin{align*}
\lambda_{k}=0  &  \Rightarrow\left\Vert f(x)\right\Vert <\delta_{k},\\
\lambda_{k}=1  &  \Rightarrow\left\Vert f(x)\right\Vert >\delta_{k+1}.
\end{align*}
Note that if $\lambda_{1}=1$, then by the continuity of $f$ at $x$, the whole sequence $\left(  z_{n}\right)  _{n\geqslant1}$ is bounded away from $x$. We may therefore assume that $\lambda_{1}=0$. If $\lambda_{k}=0$, then, using the Z-stability of $f$ at $x$, choose $\zeta_{k}\in Z_{f}$ with $\rho(\zeta_{k},x)<2^{-n}r$. If $\lambda_{k}=1-\lambda_{k-1}$, set $\zeta_{j}=\zeta_{k-1}$ for each $j\geqslant k$; then $\rho(\zeta_{k},x)<2^{-k+1}r$. It readily follows that
\[ \rho\left(  \zeta_{m},\zeta_{n}\right)  \leqslant2^{-n+1}r\ \ \ \ \left(m\geqslant n\right)  , \]
Since $\overline{B}(x,r)$ is a closed, and therefore complete, subset of $E$, the sequence $\left(  \zeta_{k}\right)  _{k\geqslant1}$ converges to a limit $\zeta_{\infty}\in\overline{B}(x,r)\subset\subset U$, and $\rho\left(\zeta_{\infty},\zeta_{n}\right)  \leqslant2^{-n+1}r$ for each $n$. By the continuity of $f$ at $\zeta_{\infty}$, we have $f(\zeta_{\infty})=0$. Our hypotheses on $Z_{f}$ now provide $N$ such that $\rho(z_{n},\zeta_{\infty})>2^{-N+2}r$ for all $n\geqslant N$. If $\lambda_{N}=1$, then $\left(z_{n}\right)  _{n\geqslant1}$ is bounded away from $x$. We may therefore assume that $\lambda_{N}=0$. Consider any $n\geqslant N$. If $\lambda_{n}=0$, then $\rho(\zeta_{n},x)<2^{-n}r$, so
\begin{align*}
\rho(z_{n},x)  &  \geqslant\rho(z_{n},\zeta_{\infty})-\rho(\zeta_{\infty},\zeta_{n})-\rho(\zeta_{n},x)\\
&  >2^{-N+2}r-2^{-N+1}r-2^{-n}r\\
&  \geqslant2^{-N+2}r-2^{-N+1}r-2^{-N}r=2^{-N}r.
\end{align*}
If $\lambda_{n}=1$, then there exists $k$ with $N<k\leqslant n$ such that $\lambda_{k}=1-\lambda_{k-1}$. In this case, $\zeta_{\infty}=\zeta_{k-1}$, so
\begin{align*}
\rho(z_{n},x)  &  \geqslant\rho(z_{n},\zeta_{\infty})-\rho(x,\zeta_{k-1})\\
&  >2^{-N+2}r-2^{-k+1}r\\
&  \geqslant2^{-N+2}r-2^{-N+1}r=2^{-N+1}r.
\end{align*}
We now see that $\rho(z_{n},x)>2^{-N+1}r$ for all $n\geqslant N$. Hence $\left(  z_{n}\right)  _{n\geqslant1}$ is eventually bounded away from each point $x$ of $U$.\footnote{{Examination of the proof of Lemma \ref{aa2} shows that the continuity of }$f${ is used}
\begin{itemize}
\item {(twice) to show that if }$\lambda_{k}=1${, and therefore }$\left\Vert f(x)\right\Vert >\delta_{k+1}${, then }$\left(  z_{n}\right)  _{n\geq1}${ is bounded away from }$x${, and}
\item {to show that }$f(\zeta_{\infty})=0${.}
\end{itemize}

{The second application can be avoided altogether by this argument. Suppose that }$f(\zeta_{\infty})\neq0${. If }$\lambda_{k}=0${\ for all }$k${, then }$\zeta_{\infty}=x${\ and }$f(x)=0${, a contradiction. If }$\lambda_{k}=1-\lambda_{k-1}${, then }$f(\zeta_{\infty})=f(\zeta_{k-1})=0${, again a contradiction. It follows from all this that if }$f(\zeta_{\infty})\neq0${, then we have neither }$\lambda_{n}=0${\ for all }$n${\ nor }$\lambda_{n}=1${\ for some }$n${. Since this is absurd, we cannot have }$f(\zeta_{\infty})\neq0${, so }$f(\zeta_{\infty})=0${.}

{In view of this, we could replace the continuity hypothesis in Lemma \ref{aa2} by the following one: for each }$x\in E${, if }$f(x)\neq 0${, then }$\rho\left(  x,Z_{f}\right)  >0${. At the same time, however, we should bear in mind that in both the recursive and the intuitionistic `models' of BISH, every function from a complete, separable metric space into a metric space is pointwise continuous everywhere; see \cite[Ch. 3, Sec 6, and Ch. 5, Corollary (2.4)]{BR}.}}
\qed

In the final result of this section we use the full \emph{anti-Specker property} (relative to one-point extensions) for a metric space $X$,
\begin{quote}
\textbf{AS}$^{X}$: If $X\cup\left\{  \omega\right\}  $ is a one-point extension of $X$, and $\left(  x_{n}\right)  _{n\geqslant1}$ is a sequence in $X\cup\left\{  \omega\right\}  $ that is eventually bounded away from each point of $X$, then there exists $N$ such that $x_{n}=\omega$ for all $n\geqslant N$.
\end{quote}
This, like \textbf{AS}$_{L}^{X}$, does not depend on the one-point extension of $X$: if it holds for some one-point extension of $X$, then it holds for them all. The statement `\textbf{AS}$^{X}$ holds for every compact metric space'\textbf{\ }is equivalent, over \textbf{BISH}, to Brouwer's fan theorem for c-bars \cite{BergerB}.

\begin{prop}
\label{2503a1}The following are equivalent over \textbf{\emph{BISH}}:
\begin{enumerate}[label=(\roman*)]
\item \textbf{\emph{AS}}$^{\left[  0,1\right]  }.$
\item \textbf{\emph{AS}}$^{2^{\mathbf{N}}}$.
\item \textbf{\emph{AS}}$^{X}$ holds for every compact metric space $X$.
\end{enumerate}
\end{prop}
\proof 
Suppose that \textbf{AS}$^{\left[  0,1\right]  }$ holds. Then (see \cite[Proposition 1]{Dien1b}), $2^{\mathbf{N}}$ has the anti-Specker property. But if $X$ is any compact metric space, then \cite[Chapter 5, Proposition (1.4)]{BR} shows that there exists a uniformly continuous mapping of $2^{\mathbf{N}}$ onto $X$; whence, by \cite[Proposition 10]{BD1}, $X$ has the anti-Specker property. Hence each of \textbf{AS}$^{\left[0,1\right]  }$ and \textbf{AS}$^{2^{\mathbf{N}}}$ is prototypical for the full anti-Specker property on compact metric spaces.\footnote{{This proposition originally appeared as \cite[Proposition 1]{Dien1b}. Our proof corrects the argument at the end of Diener's one.}}
\qed
A classical sequential compactness argument shows that under the hypotheses of Lemma \ref{aa2}, the set of zeroes of $f\,$\ in a compact set well contained in $U$ is finite. Here is \ our constructive counterpart of that result.

\begin{prop}
\label{aa3}Let $E$ be a complete metric space, $U$ an open subset of $E$, and $f$ a pointwise continuous mapping of $U$ into a normed space $Y$. Suppose that $f$ is Z-stable at each point of $U$, and that $Z_{f}$ is countably isolated, with one-one enumeration $z_{1},z_{2},\ldots$ . Let $X$ be a compact subset of $U$ with the anti-Specker property. Then there exists $N$ such that
\[ Z_{f}\cap X\subset\left\{  z_{1},\ldots,z_{N}\right\}  . \]
\end{prop}
\proof 
By Lemma \ref{aa2}, $\left(  z_{n}\right)  _{n\geqslant1}$ is eventually bounded away from each point of $X$. By the anti-Specker property for $X$, there exists $N$ such that $\rho\left(  z_{n},x\right)  \geqslant2^{-N}$ whenever $x\in X$ and $n>N$. Thus if $x\in X$ and $f(x)=0$, we must have $x=z_{n}$ for some $n\leqslant N$.
\qed

\section{Concluding remarks on \textbf{AS}$_{L}^{X}$ and positivity properties}

In order to facilitate the passage from Z-stability to uniform Z-stability in Proposition \ref{2003a1}, we replaced the positivity property, used in Proposition \ref{2303a2}, by\textbf{\ }the limited anti-Specker property. An examination of the proof of \cite[Proposition 5]{BD1} shows that \textbf{AS}$_{L}^{X}$ implies the \emph{pointwise positivity property},
\begin{quote}
\textbf{POS}$_{p}^{X}$: \ If $f$ is a pointwise continuous, positive-valued mapping on a metric space, \emph{and if }$\inf_{X}f$\emph{\ exists}, then $\inf_{X}f>0$.
\end{quote}
The following lemma will enable us to prove that the pointwise positivity property for $\left[  0,1\right]  $ implies, and hence is equivalent to, \textbf{AS}$_{L}^{\left[  0,1\right]  }$. For the lemma, we need this definition. For each $t\in\mathbf{R}$ and each $\delta>0$ we define the corresponding \emph{spike function} to be the unique uniformly continuous function $s\left(  t,\delta,.\right)  :\mathbf{R}\rightarrow\mathbf{R}$ with the following properties:
\begin{itemize}
\item[$\cdot$] $s\left(  t,\delta,t\right)  =1,$
\item[$\cdot$] $s\left(  t,\delta,x\right)  =0$ whenever $\left\vert x-t\right\vert >\delta,$ and
\item[$\cdot$] $s\left(  t,\delta,.\right)  $ is linear in each of the intervals $\left[  t-\delta,t\right]  $ and $\left[  t,t+\delta\right]  .$
\end{itemize}

\begin{lemma}
\label{puc1a}Let $\left(  x_{n}\right)  _{n\geqslant1}$ be a sequence in $\left[  0,1\right]  \cup\left\{  2\right\}  $ that is eventually bounded away from each point of $\left[  0,1\right]  .$ Let $\delta_{1}\in\left(0,\frac{1}{2}\right)  ,$ and for each $n\geqslant2$ let $0<\delta_{n}\leqslant\min\left\{  2^{-n},\delta_{n-1}\right\}  $. Let $\left(a_{n}\right)  _{n\geqslant1}$ be a sequence of real numbers, and define, for each $n$, a mapping $f_{n}:\left[  0,1\right]  \rightarrow(0,1]$ by
\[ f_{n}(x)\equiv
\begin{cases}
 s(x_{n},\delta_{n},.) & \text{if }x_{k}\in\left[  0,1\right]  \text{ for each }k\leqslant n\\
0 & \text{if }x_{k}=2\text{ for some }k\leqslant n.
\end{cases} \]
Then $f=\sum_{n=1}^{\infty}a_{n}f_{n}$ is a well-defined, pointwise continuous mapping on $\left[  0,1\right]  $.
\end{lemma}
\proof 
The proof of \cite[Lemma 4]{dsbarch} carries over \emph{mutatis mutandis}.
\qed

\begin{prop}
\label{2103a1}\textbf{\emph{POS}}$_{p}^{\left[  0,1\right]  }$ and \textbf{\emph{AS}}$_{L}^{\left[  0,1\right]  }$ are equivalent over \textbf{\emph{BISH.}}
\end{prop}
\proof 
In view of earlier remarks, it is enough to assume \textbf{POS}$_{p}^{\left[0,1\right]  }$ and derive \textbf{AS}$_{L}^{\left[  0,1\right]  }$. Consider a sequence $\left(  z_{n}\right)  _{n\geqslant1}$ in $Y\equiv\left[  0,1\right]\cup\left\{  2\right\}  $ that is eventually bounded away from each point of $\left[  0,1\right]  $. Using \cite[Lemma 5]{dsbarch}, we may take $z_{m}\neq z_{n}$ whenever $m\neq n$. Replacing the sequence $\left(z_{n}\right)  _{n\geqslant1}$ by one of its tails if necessary, we may further assume that there exists $\delta_{0}\in\left(  0,\frac{1}{2}\right)  $ such that $z_{n}>\delta_{0}$ and $\left\vert z_{n}-1\right\vert >\delta_{0}$ for all $n.$ Setting $n_{0}=1$ and arguing as in the proof of \cite[Theorem 6]{dsbarch}, construct, inductively, a sequence $\left(  \delta_{k}\right)_{k\geqslant0}$ of positive numbers and a strictly increasing sequence $\left(  n_{k}\right)  _{k\geqslant0}$ of positive integers such that the following hold for each $k\geqslant1:$
\begin{enumerate}[label=(\roman*)]
\item $\delta_{k}\leqslant\min\left\{  2^{-k},\delta_{k-1}\right\}  ;$
\item $\left\vert z_{j}-z_{k}\right\vert \geqslant2\delta_{k}$ for all $j\geqslant n_{k}.$
\end{enumerate}
Replacing $\delta_{k}$ by a smaller value if necessary, since the terms of $\left(  z_{n}\right)  _{n\geqslant1}$ are distinct we may assume that $\left\vert z_{j}-z_{k}\right\vert \geqslant2\delta_{k}$ also when $k\neq j<n_{k}$. It follows from Lemma \ref{puc1a} that
\begin{equation}
f=\sum_{k=1}^{\infty}\left(  1-2^{-k}\right)  s\left(  z_{k},\delta_{k},.\right)  \label{bb1}
\end{equation}
defines a pointwise continuous function on $\left[  0,1\right]  $. Consider any $j,k$ with $j>k$. If $\left\vert x-z_{j}\right\vert <\delta_{j}$ and $\left\vert x-z_{k}\right\vert <\delta_{k}$, then as $\delta_{j}\leqslant\delta_{k}$, we have $\left\vert z_{j}-z_{k}\right\vert <2\delta_{k}$, a contradiction from which we conclude that the supports of the terms of the series at (\ref{bb1}) are pairwise disjoint. Hence $f$ maps $\left[0,1\right]  $ into $[0,1).$ We now prove that $\sup f$ exists. Let $0<\alpha<\beta<1$, and pick a positive integer $K$ such that $\beta <1-2^{-K+1}$. Either there exists $k\leqslant K$ such that $x_{k}=2$, in which case $\sup f $ $=1-2^{-k+1}$ for the smallest such $k$, and either $\sup f>\alpha$ or $\sup f<\beta$; or else $x_{k}\in\left[  0,1\right]  $ for all $k\leqslant K$, and $f(x_{K})=1-2^{-K}>\beta>\alpha$. Since $\alpha,\beta$ are arbitrary, it follows from the constructive least-upper-bound principle \cite[Theorem 2.1.18]{BV} that $\sup f$ exists.

Now let $g\equiv1-f$, which is a pointwise continuous mapping of $\left[0,1\right]  $ into $(0,1]$ with infimum equal to $1-\sup f$. Applying \textbf{POS}$_{p}^{[0,1]}$ to $g$, we see that $\inf g>0$. Thus there exists a positive integer $\kappa$ such that $2^{-\kappa}<\inf g$. It follows that $z_{n_{k}}=x_{k}=2$ for some $k\leqslant\kappa$.
\qed

It is no surprise that each of \textbf{AS}$_{L}^{\left[  0,1\right]  }$ and \textbf{AS}$_{L}^{2^{\mathbf{N}}}$ is prototypical for the limited anti-Specker property on compact spaces.

\begin{prop}
\label{2303a4}The following are equivalent over \textbf{\emph{BISH}}:
\begin{enumerate}[label=(\roman*)]
\item \textbf{\emph{AS}}$_{L}^{\left[  0,1\right]  }$.
\item \textbf{\emph{AS}}$_{L}^{2^{\mathbf{N}}}$.
\item \textbf{\emph{AS}}$_{L}^{X}$ holds for each compact metric space $X$.
\end{enumerate}
\end{prop}
\proof 
The argument used (before Proposition \ref{aa3}) to prove that \textbf{AS}$^{2^{\mathbf{N}}}$ implies \textbf{AS}$^{X}$ for any compact $X$ trivially adapts to show that if \textbf{AS}$_{L}^{2^{\mathbf{N}}}$ holds, then every compact metric space has the limited anti-Specker property. Diener \cite{Dien13} has recently proved that \textbf{AS}$_{L}^{\left[  0,1\right]}$ implies \textbf{AS}$_{L}^{2^{\mathbf{N}}}$. Clearly, (iii) implies (i).
\qed

In turn, \textbf{POS}$_{p}^{\left[  0,1\right]  }$ and \textbf{POS}$_{p}^{2^{\mathbf{N}}}$ are prototypical in their realm:
\begin{prop}
\label{2303a5}The following are equivalent over \textbf{\emph{BISH}}:
\begin{enumerate}[label=(\roman*)]
\item \textbf{\emph{POS}}$_{p}^{\left[  0,1\right]  }$.
\item \textbf{\emph{POS}}$_{p}^{2^{\mathbf{N}}}$.
\item \textbf{\emph{POS}}$_{p}^{X}$ holds for each compact
metric space $X$.
\end{enumerate}
\end{prop}
\proof 
This follows from Propositions \ref{2103a1} and \ref{2303a4}.
\qed

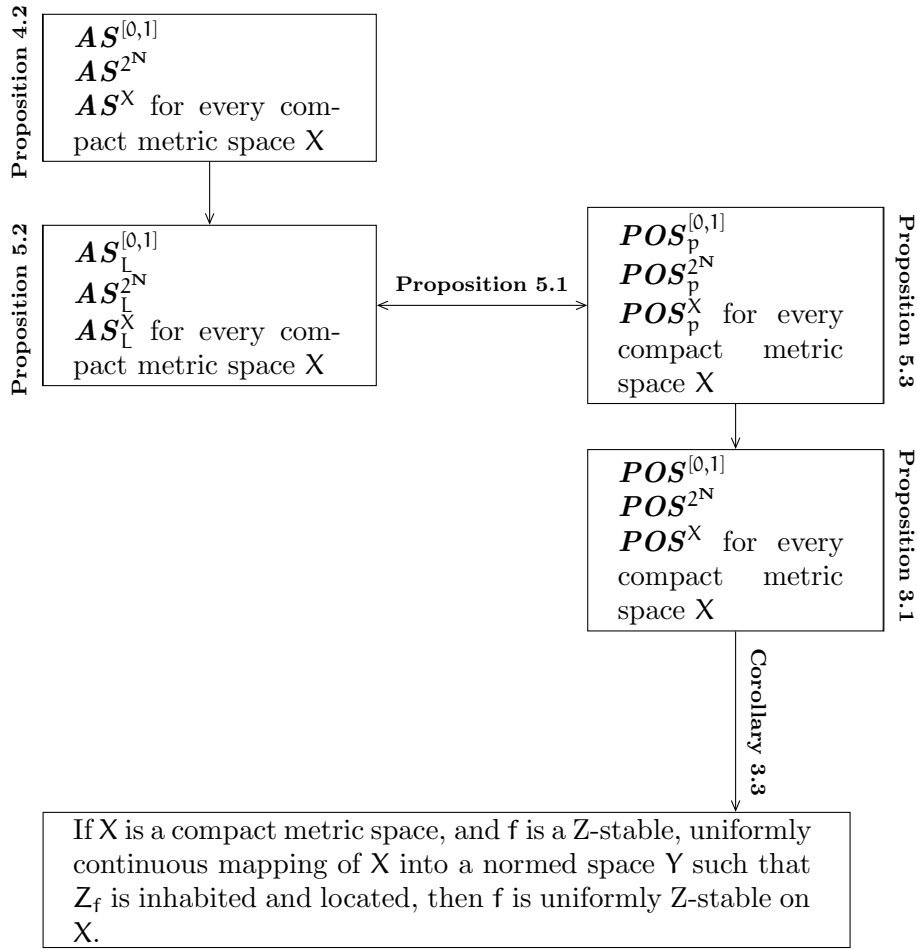
\begin{figure}[h]
  \centering
  \scalebox{1.0}{\begin{tikzpicture}
  
    \newcommand{\widestl}{$\negspeckinc[\ui]$}
    \newcommand{\widestc}{$\limaspeck$}
    \newcommand{\widestr}{$\hbcountints{\ui}$}
    
    \matrix[column sep=28mm,row sep=6mm,ampersand replacement=\&] (Layout) {
    
      \node[equivclass] (AntiSpecker) {\equivclassbox{35mm}{
        \item \textbf{\emph{AS}}$^{\left[  0,1\right]  }$
        \item \textbf{\emph{AS}}$^{2^{\mathbf{N}}}$
        \item \textbf{\emph{AS}}$^{X}$ for every compact metric space $X$
      }}; \\
      
      \node[equivclass] (LimAntiSpecker) {\equivclassbox{35mm}{
        \item \textbf{\emph{AS}}$_{L}^{\left[  0,1\right]  }$
        \item \textbf{\emph{AS}}$_{L}^{2^{\mathbf{N}}}$
        \item \textbf{\emph{AS}}$_{L}^{X}$ for every compact metric space $X$
      }}; \&
      
      \node[equivclass] (PointwisePositivity) {\equivclassbox{30mm}{
        \item \textbf{\emph{POS}}$_{p}^{\left[  0,1\right]  }$
        \item \textbf{\emph{POS}}$_{p}^{2^{\mathbf{N}}}$
        \item \textbf{\emph{POS}}$_{p}^{X}$ for every compact metric space $X$
      }}; \\
      
      \&
      \node[equivclass] (Positivity) {\equivclassbox{30mm}{
        \item \textbf{\emph{POS}}$^{\left[  0,1\right]  }$
        \item \textbf{\emph{POS}}$^{2^{\mathbf{N}}}$
        \item \textbf{\emph{POS}}$^{X}$ for every compact metric space $X$
      }}; \\
    };
    
    
    \node[equivclass,anchor=north west] (ZSToUZS) at ($
    (LimAntiSpecker.west |- Positivity.south) - (0, 24 mm) $) {\equivclassbox{98mm}{
      \item If $X$ is a compact metric space, and $f$ is a Z-stable, uniformly continuous mapping of $X$ into a normed space $Y$ such that $Z_{f}$ is inhabited and located, then $f$ is uniformly Z-stable on $X$.
    }};
    
    \draw[impl] (AntiSpecker) -- (LimAntiSpecker);
    \draw[biimpl] (LimAntiSpecker) -- node[citation,above] {Proposition \ref{2103a1}} (PointwisePositivity);
    \draw[impl] (PointwisePositivity) -- (Positivity);
    \Hyp{(Positivity.south)}{(Positivity.south |- ZSToUZS.north)}{Corollary \ref{2303a3}} ;
\HyP{(AntiSpecker.south west)}{(AntiSpecker.north west)}{\!Proposition \ref{2503a1}\,};
\HyP{(LimAntiSpecker.south west)}{(LimAntiSpecker.north west)}{\!Proposition \ref{2303a4}\,};
\HyP{(PointwisePositivity.north east)}{(PointwisePositivity.south east)}{\!Proposition \ref{2303a5}\,};
\HyP{(Positivity.north east)}{(Positivity.south east)}{\!Proposition \ref{2303a1}\,};

  \end{tikzpicture}}
  \bigskip
  \caption{{Relationships between anti-Specker and positivity equivalence classes, relative to \textbf{BISH}.}}
  \label{fig:ASPosRelationships}
\end{figure}

\textbf{Figure \ref{fig:ASPosRelationships}} summarises the relationships between the anti-Specker and positivity properties that we have established here: clearly, \textbf{POS}$_{p}^{\left[  0,1\right]  }$, and therefore \textbf{AS}$_{L}^{\left[  0,1\right]  }$, implies \textbf{POS}$^{\left[0,1\right]  }$. Does \textbf{POS}$^{\left[  0,1\right]  }$ imply \textbf{AS}$_{L}^{\left[  0,1\right]  },$ and therefore \textbf{AS}$_{L}^{X}$ for every compact metric space $X$? We do not know the answer; but if it turns out to be `yes', then in view of the Julian-Richman theorem \cite{JR}, we will have found the exact fan-theoretic equivalent, relative to \textbf{BISH}, of \textbf{AS}$_{L}^{[0,1]}$: namely, \textbf{FT}$_{D}$.
\let\c2303a3={Corollary \ref{2303a3}}

\section*{Acknowledgements}
The authors thank the University of Canterbury for supporting Dent by a doctoral scholarship during the writing of this paper, and the Royal Society of New Zealand for partial support of McKubre-Jordens by means of a Marsden Fund grant.

\end{document}